\documentclass[12pt,reqno]{amsart}
\usepackage[english]{babel}
\usepackage{amsmath,amsfonts,amssymb,amsthm,amscd,latexsym}
\usepackage[T1]{fontenc}
\usepackage[utf8]{inputenc}

\usepackage{tikz}
\usetikzlibrary{arrows}
\usepackage{color}
\usepackage{cite}
\usepackage{multirow}
\usepackage{cite}

\newtheorem{theorem}{Theorem}[section]

\newtheorem{lemma}[theorem]{Lemma}

\newtheorem{definition}[theorem]{Definition}
\newtheorem{exam}[theorem]{Example}

\numberwithin{equation}{section}
\DeclareFontFamily{OT1}{pzc}{}
\DeclareFontShape{OT1}{pzc}{m}{it}{<-> s * [1.10] pzcmi7t}{}
\DeclareMathAlphabet{\mathpzc}{OT1}{pzc}{m}{it}

\def\ms{\mathfrak{s}}

\def\mh{\mathfrak{h}}
\def\so{\mathfrak{so}}
\def\sl{\mathfrak{sl}}
\def\gl{\mathfrak{gl}}

\def\C{\Bbb C}

\begin{document}


\title[Local derivation on the Schr{\"o}dinger Lie algebra  in $(n+1)$-dimensional space-time]{Local derivation on the Schr{\"o}dinger Lie algebra  in $(n+1)$-dimensional space-time}

\author{ Alauadinov A.K.,\ Yusupov B.B.}
\address[Amir Alauadinov]{
Ch.Abdirov 1, Department of Mathematics, Karakalpak State University, Nukus 230113, Uzbekistan
\newline
and\newline	
V.I.Romanovskiy Institute of Mathematics, Uzbekistan Academy of Sciences, Univesity Street, 9, Olmazor district, Tashkent, 100174, Uzbekistan		
}
	\email{amir\_t85@mail.ru}

	\address[Bakhtiyor Yusupov]{	
		V.I.Romanovskiy Institute of Mathematics, Uzbekistan Academy of Sciences, Univesity Street, 9, Olmazor district, Tashkent, 100174, Uzbekistan
\newline
and\newline
Department of Physics and Mathematics, Urgench State University, H. Alimdjan street, 14, Urgench
220100, Uzbekistan
}
\email{baxtiyor\_yusupov\_93@mail.ru\ b.yusupov@mathinst.uz}

\date{}
\maketitle

\begin{abstract} This paper studies local derivations on the Schr{\"o}dinger algebra $\ms_n$ in $(n+1)$-dimensional space-time of Schr{\"o}dinger Lie groups for any integer $n$. The purpose of this work is to  prove that every local derivation on $\ms_n$  is a derivation.

\\


\end{abstract}

\maketitle \thispagestyle{empty}

\section{Introduction}\label{sec:intro}

Local derivations are useful tools in studying the structure of rings and algebras, where there are still many related unsolved problems. R.V.Kadison, D.R.Larson and A.R.Sourour first introduced the notion of local derivations on algebras in their remarkable paper \cite{Kad, Lar}. Since then many researches have been studying local derivations of different types of algeras (e.g. see \cite{AyuKud, AyuKhu, KudOmiKur, AlaYus1}). In \cite{AyuKud} the authors proved that every local derivation  on a finite-dimensional semisimple Lie algebra $\mathcal{L}$ over an algebraically closed field of characteristic zero is a derivation. In \cite{AyuKhu}  local derivations of solvable Lie algebras are investigated and it is shown that in the class of solvable Lie algebras, there exist algebras which admit local derivations which are not ordinary derivation and also algebras for which every local derivation is a derivation. Moreover, it is proved that every local derivation on a finite-dimensional solvable Lie algebra with model nilradical and maximal dimension of complementary space is a derivation. In \cite{KudOmiKur}, the authors proved that every  local derivations on solvable Lie algebras whose nilradical has maximal rank is a derivation. In \cite{AlaYus1}, the authors proved  that every local derivation on the conformal Galilei algebra is a derivation. 

We note that aforementioned algerbas are finite-dimensional algebras. In the infinite-dimensional case, 
the authors of \cite{Kaim,AyuKudYus,Yao} proved that every local derivation on some class of the locally simple Lie algebras, generalized Witt algebras, Witt algebras, the Witt algebras over a field of prime characteristic is a derivation.

The Schrödinger Lie group describes symmetries of the free particle Schrödinger equation in \cite{Per}. The Lie algebra $\mathfrak{s}_n$ in $(n+1)$-dimensional space-time of the Schrödinger Lie group is called the Schrödinger algebra, see \cite{DDM1}. The Schrödinger algebra $\mathfrak{s}_n$ is a non-semisimple Lie algebra and plays an important role in mathematical physics. Recently there was a series of papers on studying the structure and representation theory of the Schrödinger algebra $\mathfrak{s}_1$ in the case of $(1+1)$-dimensional space-time, see \cite{AD, YT, ZCH}.	
	 	
In this paper, we will generalize our previous result for all $n>2$  where we proved that every local derivation on the Schrödinger algebra $\mathfrak{s}_n$ in $(n+1)$-dimensional
space-time is a derivation when $n = 1,2$ \cite{AlaYus2}. Hence, the result will be the same for all $n\in\mathbb{N}$.

\section{Preliminaries}

In this section, we first recall the definition of $\ms_n$ from \cite{DDM1} in a different form. We know that  the general linear Lie algebra $\mathfrak{gl}_{2n}$  has  the natural representation on $\C^{2n}$ by left matrix  multiplication. Let $\{e_1, e_2, \cdots, e_{2n}\}$ be  the standard basis of $\C^{2n}$.

The Heisenberg Lie algebra $\mh_{n}=\C^{2n}\oplus \C z$ is the Lie algebra with
Lie bracket   given by
$$[e_i, e_{n+i}]=z, \qquad [z, \mh_n]=0,\quad 1\leq i\leq n.$$

Recall that the Schr{\"o}dinger Lie   algebra $\ms_n$ is  the  semidirect product Lie algebra
$$\ms_n=(\sl_2\oplus \so_n )\ltimes \mh_n,$$
where $\sl_2$ is embedded in $\gl_{2n}$
by the mapping
$$\left(
    \begin{array}{cc}
      a & b \\
      c & -a \\
    \end{array}
  \right)\mapsto \left(
    \begin{array}{cc}
      a I_n& bI_n \\
      cI_n & -aI_n \\
    \end{array}
  \right)
$$
and $\so_n$ is embedded in $\gl_{2n}$ by
$$A\in \so_n\mapsto \left(
    \begin{array}{cc}
     A& 0 \\
      0 & A \\
    \end{array}
  \right). $$

  Here $I_n$ is the $n\times n$ identity matrix, $\sl_2\oplus \so_n$ acts
   on $\mh_n$ by matrix multiplication, and $[z,\ms_n]=0$ .

Next, we will introduce a basis of $\ms_n$. Let $$\aligned h&=\left(
    \begin{array}{cc}
       I_n& 0 \\
       0 & -I_n \\
    \end{array}
  \right),
 e=\left(
    \begin{array}{cc}
      0& I_n \\
      0& 0\\
    \end{array}
  \right),
  f=\left(
    \begin{array}{cc}
      0& 0\\
      I_n & 0\\
    \end{array}
  \right),\\
s_{ij}&=\left(
    \begin{array}{cc}
     e_{ij}-e_{ji}& 0 \\
      0 &  e_{ij}-e_{ji} \\
    \end{array}  \right),\quad 1\leq i<j\leq n,\\
    u_k&=e_k, v_k=e_{n+k},\quad 1\leq k\leq n, \endaligned $$where $e_{i,j}(1\leq i,j\leq n)$ the $n\times n$ matrix with zeros everywhere except a $1$ on position $(i,j)$.

Note that $s_{ij}=-s_{ji}$. We can check that the  non-trivial commutation relations of $\ms_n$ are

The Schr\"{o}dinger algebra $\mathfrak{s}_{n}$ is a Lie algebra with a $\mathbb{C}$-basis
$$
\{e,f,h,z,u_{i},v_{i},s_{jk}(=-s_{kj})\mid1\leqslant i\leqslant n,1\leqslant j<k\leqslant n\}
$$
equipped with the following non-trivial commutation relations
\begin{equation}\label{Bracket}
\aligned & [h,e]=2e, [h,f]=-2f, [e,f]=h,\\
& [u_i, v_i]=z, [h,u_i]=u_i, [h,v_i]=-v_i,\\
& [e,v_i]=u_i, [f,u_i]=v_i,\\
&[s_{kl},  u_i]=\delta_{li} u_k-\delta_{ki}u_l, [s_{kl},  v_i]=\delta_{li} v_k-\delta_{ki}v_l, \\
& [s_{ij},s_{kl}]=\delta_{kj}s_{il}+\delta_{il}s_{jk}+\delta_{lj}s_{ki}+\delta_{ki}s_{lj}.
\endaligned
\end{equation}
where $\delta_{ij}$ is the Kronecker Delta defined as $1$ for $i=j$ and as $0$ otherwise.

 The Schr{\"o}dinger algebra $\mathfrak{s}_{n}$ is a finite-dimensional, non-semisimple and non-solvable Lie algebra, and it is the semidirect product Lie algebra
$$\mathfrak{s}_{n}=(\mathfrak{sl}_{2}\oplus\mathfrak{so}_{n})\ltimes\mathfrak{h}_{n}$$
where $\mathfrak{sl}_{2}={\rm Span}_{\mathbb{C}}\{e,f,h\}$ is the 3-dimensional simple Lie algebra, $\mathfrak{so}_{n}={\rm Span}_{\mathbb{C}}\{s_{kl}\mid1\leqslant k<l\leqslant n\}$ is the orthogonal Lie algebra and $\mathfrak{h}_{n}={\rm Span}_{\mathbb{C}}\{z,u_{i},v_{i}\mid1\leqslant i\leqslant n\}$ is the Heisenberg Lie algebra.  

A derivation on a  Lie algebra $\mathcal{L}$ is a linear map  $D: \mathcal{L} \rightarrow \mathcal{L}$  which satisfies the Leibniz rule:
\begin{equation}\label{Defder}
D([x,y]) = [D(x), y] + [x, D(y)], \quad \text{for any} \quad x,y \in \mathcal{L}.
\end{equation}
The set of all derivations of $\mathcal{L}$ is denoted by $\mathrm{Der}(\mathcal{L})$ and with respect to the commutation operation is a Lie algebra.

For any element $y\in \mathcal{L}$ the operator of
right multiplication $\mathrm{ad}_y: \mathcal{L} \to \mathcal{L}$, defined as $\mathrm{ad}_y(x)=[y,x]$ is a derivation, and derivations of this form are called inner derivation.
The set of all inner derivations of $\mathcal{L},$ denoted by $\mathrm{Inn}(\mathcal{L}),$ is an ideal in $\mathrm{Der}(\mathcal{L}).$

\begin{definition}
A linear operator $\Delta$ is called a local derivation if for any $x \in \mathcal{L},$ there exists a derivation $D_x: \mathcal{L} \rightarrow \mathcal{L}$ (depending on $x$) such that
$\Delta(x) = D_x(x).$ The set of all local derivations on $\mathcal{L}$ we denote by $\mathrm{LocDer}(\mathcal{L}).$
\end{definition}

We use the following definition given in \cite{QT}.

\begin{definition}\label{def1} 
The following derivations are outer derivations of $\ms_n$.
\begin{itemize}
  \item When $n\geq2$, the derivation $\sigma:\ms_n\rightarrow \ms_n$ is given by
$$\sigma(f)=\sigma(e)=\sigma(h)=\sigma(s_{kl})=0,\sigma(z)=z,\sigma(u_i)=\frac{1}{2}u_i,\sigma(v_i)=\frac{1}{2}v_i,$$      for all $1\leq i\leq n,\ 1\leq k<l\leq n.$
  \item When $n=2$, the derivation $\tau:\ms_2\rightarrow \ms_2$ is given by
  $$\tau(f)=\tau(e)=\tau(h)=\tau(z)=\tau(u_i)=\tau(v_i)=0,\tau(s_{12})=z,\ i=1,2.$$
  \item When $n=1$, the derivation $\sigma_1:\ms_1\rightarrow \ms_1$ is given by
  $$\sigma_1(f)=\sigma_1(e)=\sigma_1(h)=0,\sigma_1(z)=z,\sigma_1(u_1)=\frac{1}{2}u_1,\sigma_1(v_1)=\frac{1}{2}v_1.$$
\end{itemize}
\end{definition}

The following theorem is proved in \cite{QT}.

\begin{theorem}\label{th23} The derivations of the Schr{\"o}dinger algebra $\ms_n$ are given by
$$\mathrm{Der}(\ms_n)=\left\{
  \begin{array}{ll}
    \mathrm{Inn}(\ms_1)\oplus\mathbb{C}\sigma_1, & n=1; \\
    \mathrm{Inn}(\ms_2)\oplus\mathbb{C}\sigma\oplus\mathbb{C}\tau, & n=2; \\
    \mathrm{Inn}(\ms_n)\oplus\mathbb{C}\sigma, & n>2.
  \end{array}
\right.$$

where $\sigma_1,\tau,\sigma$ are given by Definition \ref{def1}.
\end{theorem}

The following example plays a vital role in our calculations:
\begin{exam}
The matrix form of the derivations of $\mathfrak{s}_3$ is as follows:   
\end{exam}

\[{\small
\left(
              \begin{array}{ccccccccccccc}
            2a_h & 0     &-2a_e & 0 &  0   & 0   & 0           & 0        & 0         & 0   & 0          & 0        & 0 \\
            0    &-2a_h  & 2b_f & 0 &  0   & 0   & 0           & 0        & 0         & 0   & 0          & 0        & 0 \\
            -a_f & a_e   & 0 & 0 &  0   & 0   & 0           & 0        & 0         & 0   & 0          & 0        & 0 \\
            0    & 0     & 0 & \lambda &  -a_{v_1}   & -a_{v_2}   & -a_{v_3}           & a_{u_1}        & a_{u_2}         & a_{u_3}   & 0          & 0        & 0 \\
        -a_{v_1} & 0     & -a_{u_1} & 0 &  a_h+\frac{\lambda}{2}   & a_{s_{11}}   &  a_{s_{12}}          & a_e        & 0         & 0   & -a_{v_2}          & -a_{v_3}        & 0 \\
        -a_{v_2} & 0     & -a_{u_2} & 0 &  -a_{s_{11}}    & a_h+\frac{\lambda}{2}   & a_{s_{13}}           & 0        & a_e        & 0   & a_{v_1}          & 0        & -a_{v_3} \\
        -a_{v_3} & 0     & -a_{u_3} & 0 &  -a_{s_{12}}   & -a_{s_{13}}   & a_h+\frac{\lambda}{2}           & 0        & 0         & a_e   & 0          & a_{v_1}        & a_{v_2} \\
            0    &-a_{u_1} & a_{v_1} & 0 &  a_f   & 0   & 0           & \frac{\lambda}{2}-a_h        & a_{s_{11}}         & a_{s_{12}}   & -a_{u_2}          & -a_{u_3}        & 0 \\
            0    & -a_{u_2} & a_{v_2} & 0 &  0   & a_f   & 0   & -a_{s_{11}}        &  \frac{\lambda}{2}-a_h                & a_{s_{13}}   & a_{u_1}          & 0        & -a_{u_3} \\
            0    & -a_{u_3} & a_{v_3} & 0 &  0   & 0   & a_f           & -a_{s_{12}}        & -a_{s_{13}}         & \frac{\lambda}{2}-a_h   & 0          & a_{u_1}        & a_{u_2} \\
            0    & 0 & 0 & 0 &  0   & 0   & 0           & 0        & 0         & 0   & 0          & a_{s_{13}}        & -a_{s_{12}} \\
            0    & 0 & 0 & 0 &  0   & 0   & 0           & 0        & 0         & 0   & -a_{s_{13}}          & 0        & a_{s_{11}} \\
            0    & 0 & 0 & 0 &  0   & 0   & 0           & 0        & 0         & 0   & a_{s_{12}}         & -a_{s_{11}}        & 0 \\
            \end{array}
                   \right)}
\]

Let $D:\mathfrak{s}_3\to \mathfrak{s}_3$ is a derivation.
From \eqref{Bracket} we obtain

\begin{equation}\label{Der}
D=\left( \begin{matrix}
   D_{11} & 0 & 0  \\
   D_{21} & D_{22} & D_{23}  \\
   0 & 0 & D_{33}  \\
\end{matrix} \right)
\end{equation}
where
\[
D_{11}=\left(
              \begin{array}{ccc}
2a_h & 0 & -2a_e\\
0 & -2a_h & 2b_f\\
-a_f & a_e & 0\\
\end{array}
                   \right), \quad 
D_{33}=\left(
              \begin{array}{ccc}
0 & a_{s_{13}} & -a_{s_{12}}\\
-a_{s_{13}} & 0 & a_{s_{11}}\\
a_{s_{12}} & -a_{s_{11}} & 0\\
\end{array}
                   \right),
\]
\[
D_{21}=\left(
              \begin{array}{ccc}
    0    & 0     & 0  \\
        -a_{v_1} & 0     & -a_{u_1} \\
        -a_{v_2} & 0     & -a_{u_2} \\
        -a_{v_3} & 0     & -a_{u_3}  \\
            0    &-a_{u_1} & a_{v_1}  \\
            0    & -a_{u_2} & a_{v_2} \\
            0    & -a_{u_3} & a_{v_3}  \\
            \end{array}
                   \right),
\quad
D_{23}=
\left(
              \begin{array}{ccc}
           0          & 0        & 0 \\
         -a_{v_2}          & -a_{v_3}        & 0 \\
        a_{v_1}          & 0        & -a_{v_3} \\
         0          & a_{v_1}        & a_{v_2} \\
         -a_{u_2}          & -a_{u_3}        & 0 \\
             a_{u_1}          & 0        & -a_{u_3} \\
             0          & a_{u_1}        & a_{u_2} \\
           \end{array}
                   \right),
\]

\[D_{22}=
\left(
              \begin{array}{ccccccc}
             \lambda &  -a_{v_1}   & -a_{v_2}   & -a_{v_3}           & a_{u_1}        & a_{u_2}         & a_{u_3}    \\
        0 &  a_h+\frac{\lambda}{2}   & a_{s_{11}}   &  a_{s_{12}}          & a_e        & 0         & 0    \\
         0         &-a_{s_{11}}    & a_h+\frac{\lambda}{2}   & a_{s_{13}}           & 0        & a_e        & 0    \\
      0         & -a_{s_{12}}   & -a_{s_{13}}   & a_h+\frac{\lambda}{2}           & 0        & 0         & a_e    \\
           0 &  a_f   & 0   & 0           & \frac{\lambda}{2}-a_h        & a_{s_{11}}         & a_{s_{12}}    \\
           0 &  0   & a_f   & 0   & -a_{s_{11}}        &  \frac{\lambda}{2}-a_h                & a_{s_{13}}  \\
           0 &  0   & 0   & a_f           & -a_{s_{12}}        & -a_{s_{13}}         & \frac{\lambda}{2}-a_h    \\
          
            \end{array}
                   \right).
\]
\section{Main results}

In this section, we will prove that every local derivation on the Schr{\"o}dinge algebra $\ms_n$ is a derivation.

\begin{theorem}\label{thm21}
Every local derivation on the Schr{\"o}dinge algebra $\ms_n,\ n \geq3$ is a derivation.
\end{theorem}

To obtain this result, we first have to prove a few lemmas.

\begin{lemma}\label{Lem1}
Every local derivation on $\mathfrak{sl}_2\oplus\mathfrak{so}_n$ is a derivation.
\end{lemma}
\begin{proof}
$\mathfrak{sl}_2\oplus\mathfrak{so}_n$ is a semi-simple algebra. In work \cite{AyuKud} it is proved that any local derivation on finite-dimensional semi-simple algebra is a derivation. Then every local derivation on $\mathfrak{sl}_2\oplus\mathfrak{so}_n$ is a derivation.
\end{proof}

Let $D:\mathfrak{s}_n\to \mathfrak{s}_n$ is a derivation.
From \eqref{Bracket} we obtain

\begin{equation}\label{Der}
D=\left( \begin{matrix}
   D_{11} & 0 & 0  \\
   D_{21} & D_{22} & D_{23}  \\
   0 & 0 & D_{33}  \\
\end{matrix} \right).
\end{equation}


From \eqref{Der}, every local derivation $\Delta$ on $\mathfrak{s}_n$ is of the form:
\begin{equation}\label{LDer}
\Delta=\left( \begin{matrix}
   \Delta_{11} & 0 & 0  \\
   \Delta_{21} & \Delta_{22} & \Delta_{23}  \\
   0 & 0 & \Delta_{33}  \\
\end{matrix} \right).
\end{equation}

From Lemma \ref{Lem1}, the elements $\Delta_{11}$ and $\Delta_{33}$ of the matrix of the local derivation $\Delta$ corresponds with the elements $D_{11}$ and $D_{33}$ of the matrix of the derivation $D$ i.e. $\Delta_{11}=D_{11}$ and $\Delta_{33}=D_{33}$.

Now we consider the local derivation of the form
\begin{equation}\label{LDer1}
    \Delta'=\Delta-D.
\end{equation} 
Using equalities $\Delta_{11}=D_{11}$ and  $\Delta_{33}=D_{33}$ we can derive
\begin{equation}\label{LDer}
\Delta'=\left( \begin{matrix}
   0 & 0 & 0  \\
   \Delta'_{21} & \Delta'_{22} & \Delta'_{23}  \\
   0 & 0 & 0  \\
\end{matrix} \right)
\end{equation}

For any $x\in\mathfrak{s}_n$ there exist
\[
a=a_{f}f+a_{h}h+a_{e}e+a_{z}z+\sum\limits_{i=1} ^na_{u_i}u_i+\sum\limits_{i=1}^na_{v_i}v_i+\sum\limits_{1\leq k<l\leq n}a_{s_{k,l} }s_{k,l}
\]
and $\lambda\in \mathbb{C}$ such that, by using Theorem \ref{th23} we can write
$$\Delta'(x)=[a,x]+\lambda\sigma(x),$$
where $a_e,a_f,a_h,a_z,a_{u_i},a_{v_i},a_{s_{k,l}},\lambda $ are complex numbers depending on $x.$

Considering \eqref{LDer1} instead of $x$ we put $h$ and $z$ we get:
\begin{equation}\begin{split}\label{3,5}
\Delta'(h)&=-\sum\limits_{i=1}^na^{(h)}_{u_i}u_i+\sum\limits_{i=1}^na^{(h)}_{v_i}v_i,\\
\Delta'(z)&=\lambda^{(z)}z.\\
\end{split}
\end{equation}
Let $x_0=\sum\limits_{i=1}^na^{(h)}_{u_i}u_i+\sum\limits_{i=1}^na^{(h)}_{v_i}v_i.$
Consider the following statement

\begin{equation}\begin{split}\label{L9}
\Delta''=\Delta'-\mathrm{ad}(x_0)-\lambda^{(h)} \sigma.
\end{split}
\end{equation}
Then $\Delta''$ is a local derivation. By direct verification we have
\begin{equation}\label{L37}
\Delta''=\left( \begin{matrix}
   0 & 0 & 0  \\
   \Delta''_{21} & \Delta''_{22} & \Delta''_{23}  \\
   0 & 0 & 0  \\
\end{matrix} \right)
\end{equation}
and
\begin{equation}\begin{split}\label{L11}
\Delta''(h)=\Delta''(z)=0.
\end{split}
\end{equation}

Considering (\ref{L37}), we find the values of the operator $\Delta''$ in the basis elements:
\begin{equation}\label{asos}
\begin{aligned}
\Delta''(f)&=-\sum\limits_{i=1}^na^{(f)}_{u_i}v_i,\\
\Delta''(e)&=-\sum\limits_{i=1}^na^{(e)}_{v_i}u_i,\\
\Delta''(u_i)&=a^{(u_i)}_{f}v_i+\left(a^{(u_i)}_{h}+\frac{\lambda^{(u_i)}}{2}\right)u_i-a^{(u_i)}_{v_i}z+\sum\limits_{ 1\leq k<i}a^{(u_i)}_{s_{k,i}}u_k-\sum\limits_{ i<l\leq n}a^{(u_i)}_{s_{i,l}}u_l,\\
\Delta''(v_i)&=\left(\frac{\lambda^{(v_i)}}{2}-a^{(v_i)}_{h}\right)v_i+a^{(v_i)}_{e}u_i+a^{(v_i)}_{u_i}z+\sum\limits_{ 1\leq k<i}a^{(v_i)}_{s_{k,i}}v_k-\sum\limits_{ i<l\leq n}a^{(v_i)}_{s_{i,l}}v_l,\\
\Delta''(s_{k,l})&=-a^{(s_{k,l})}_{u_l}u_k+a^{(s_{k,l})}_{u_k}u_l-a^{(s_{k,l})}_{v_l}v_k+a^{(s_{k,l})}_{v_k}v_l.
\end{aligned}
\end{equation}
We take an element
$b=b_{e}e+b_{f}f+b_{h}h+b_{z}z+\sum\limits_{i=1} ^nb_{u_i}u_i+\sum\limits_{i=1}^nb_{v_i}v_i+\sum\limits_{1\leq k<l\leq n}b_{s_{k,l} }s_{k,l}$ and $\mu\in\mathbb{C}$ for every $x\in\mathfrak{s}_n$, where $b_e,b_f,b_h,b_z,b_{u_i},b_{v_i},b_{s_{k,l}},\mu $ are complex numbers depending on $x.$

\begin{lemma}\label{Lem3}
Coefficients $a^{(u_i)}_{f}$ and $a^{(v_i)}_{e}$ $(1\leq i\leq n)$ in the formula (\ref{asos}) are equal to zero.

\end{lemma}

\begin{proof}
From the equalities, take an element $x=e+u_i$, for any fixed $1\leq i\leq n,$ we consider, using the rate of local derivation we calculate the following
\begin{equation*}\begin{split}
\Delta''(x)&=\Delta''(e+u_i)=[b,e+u_i]+\mu\sigma(e+u_i)=\\
             &=\big[b_{f}f+b_{h}h+b_{e}e+b_{z}z+
             \sum\limits_{j=1}^nb_{u_j}u_j+
             \sum\limits_{j=1}^nb_{v_j}v_j+\sum\limits_{1\leq k<l\leq n}b_{s_{k,l}}s_{k,l},\ e+u_i\big]+\\
             &+\mu\sigma(e+u_i)=-b_{f}h+b_{f}v_i+\ast e+\sum\limits_{j=1}^n\ast u_j+\ast z.\\
\end{split}
\end{equation*}

On the other hand, based on (\ref{asos}), we calculate the following equality:  
\begin{equation*}\begin{split}
\Delta''(x)&=\Delta''(e+u_i)=\Delta''(e)+\Delta''(u_i)=a^{(u_i)}_{f}v_i+\sum\limits_{j=1}^n\ast u_j+\ast z.
\end{split}
\end{equation*}

Comparing the coefficients at the basis elements $h$ and $v_i,$ we get $b_f=0,\ b_f=a^{(u_i)}_{f}$ which implies
$$a^{(u_i)}_{f}=0.$$

From the equalities, take an element $x=f+v_i$, for fixed $1\leq i\leq n,$ similarly, considering
\begin{equation*}\begin{split}
\Delta''(x)&=\Delta''(f+v_i)=[b,f+v_i]+\mu\sigma(f+v_i)=b_{e}h+b_{e}u_i+\ast f+\sum\limits_{j=1}^n\ast v_j+\ast z.
\end{split}
\end{equation*}

On the other hand,
\begin{equation*}\begin{split}
\Delta''(x)&=\Delta''(f+v_i)=\Delta''(f)+\Delta''(v_i)=a^{(v_i)}_{e}u_i+\sum\limits_{j=1}^n\ast v_j+\ast z.
\end{split}
\end{equation*}

Comparing the coefficients at the basis elements $h$ and $u_i,$ we get $b_e=0,\ b_e=a^{(v_i)}_{e},$ which implies
$$a^{(v_i)}_{e}=0.$$
\end{proof}

\begin{lemma}\label{Lem4}
 Coefficients $a^{(u_i)}_{v_i}$ and $a^{(v_i)}_{u_i}, 1\leq i\leq n$ in the formula (\ref{asos}) are equal to zero.
\end{lemma}

\begin{proof} From the equalities, take an element $x=h+u_i$ for fixed $1\leq i\leq n,$ we consider,
\begin{equation*}\begin{split}
\Delta''(x)&=\Delta''(h+u_i)=[b,h+u_i]+\mu\sigma(h+u_i)=\\
               &=2b_{f}f+\sum\limits_{j=1}^nb_{v_j}v_j+b_{f}v_i-b_{v_i}z+\ast e+\sum\limits_{j=1}^n\ast u_j.\
\end{split}
\end{equation*}

On the other hand,
\begin{equation*}\begin{split}
\Delta''(x)&=\Delta''(h+u_i)=\Delta''(h)+\Delta''(u_i)=-a^{(u_i)}_{v_i}z+\sum\limits_{j=1}^n\ast u_j.
\end{split}
\end{equation*}

Comparing the coefficients at the basis elements $f,\ z$ and $v_i,$ we get $b_{f}=b_{v_i}=0,\ b_{v_i}=a^{(u_i)}_{v_i},$
 which implies $$a^{(u_i)}_{v_i}=0.$$

From the equalities, take an element $x=h+v_i,$ for fixed $1\leq i\leq n,$ similarly, considering
\begin{equation*}\begin{split}
\Delta''(x)&=\Delta''(h+v_i)=[b,h+v_i]+\mu \sigma(h+v_i)=\\
               &=-2b_{e}e-\sum\limits_{j=1}^nb_{u_j}u_j+b_{e}u_i+b_{u_i}z+\ast f+\sum\limits_{j=1}^n\ast v_j.
\end{split}
\end{equation*}

On the other hand,
\begin{equation*}\begin{split}
\Delta''(x)&=\Delta''(h+v_i)=\Delta''(h)+\Delta''(v_i)=a_{e,u_i}u_i+a^{(v_i)}_{u_i}z+\sum\limits_{ j=1}^n \ast v_j.
\end{split}
\end{equation*}

Comparing the coefficients at the basis elements $e,\ z$ and $u_i,$ we get $b_{e}=b_{u_i}=0,\ b_{u_i}=a^{(v_i)}_{u_i},$
 which implies $$a^{(v_i)}_{u_i}=0.$$
\end{proof}

\begin{lemma}\label{Lem5}
$\Delta''(f)=0$ and $a_{h}^{(v_i)}-\frac{\lambda^{(v_i)}}{2}=0$ in the formula (\ref{asos}).\\
\end{lemma}

\begin{proof}
From the equalities, take an element $x=f-\frac{1}{2}z+y_i$
\begin{equation}\begin{split}\label{f1}
\Delta''(x)&=\left[b,f-\frac{1}{2}z+y_i\right]+\mu \sigma(f-\frac{1}{2}z+y_i)=-2b_{h}f+b_{e}h-\sum\limits_{j=1}^nb_{u_j}v_j-\frac{\mu}{2}z-\\
             &-b_{h}v_i+b_{e}u_i+b_{u_i}z+\sum\limits_{ 1\leq k<i}b_{s_{k,i}}v_k-\sum\limits_{ i<l\leq n}b_{s_{i,l}}v_l+\frac{\mu}{2}v_i.\\
\end{split}
\end{equation}

On the other hand,
\begin{equation}\begin{split}\label{f2}
\Delta''(x)&=\Delta''(f)-\Delta''(\frac{z}{2})+\Delta''(v_i)=-\sum\limits_{i=1}^na_{u_i}^{(f)}v_i-a_{h}^{(v_i)}v_i+\\
           &+\sum\limits_{ 1\leq k<i}a_{s_{k,i}}^{(v_i)}v_k-\sum\limits_{ i<l\leq n}a_{s_{i,l}}^{(v_i)}v_l+\frac{\lambda^{(v_i)}}{2}v_i.
\end{split}
\end{equation}

Comparing the coefficients at the basis elements $f,\ z$ and $v_i$, (\ref{f1}) and (\ref{f2}) we obtain that

\begin{equation}\label{f3}
a_{u_i}^{(f)}=-a_{h}^{(v_i)}+\frac{\lambda^{(v_i)}}{2}.
\end{equation}

From the equalities, take an element $x=f-\frac{1}{2}z-v_i$
\begin{equation}\begin{split}\label{f4}
\Delta''(x)&=\left[b,f-\frac{1}{2}z-v_i\right]+\mu\sigma(f-\frac{1}{2}z-v_i)=-2b_{h}f+b_{e}h-\sum\limits_{j=1}^nb_{u_j}v_j-\\
             &-\frac{\mu}{2}z+b_{h}v_i-b_{e}u_i-b_{u_i}z-\sum\limits_{ 1\leq k<i}b_{s_{k,i}}v_k+\sum\limits_{ i<l\leq n}b_{s_{i,l}}v_l-\frac{\mu}{2}v_i.\\
\end{split}
\end{equation}

On the other hand,
\begin{equation}\begin{split}\label{f5}
\Delta''(x)&=\Delta''(f)-\Delta''(\frac{z}{2})-\Delta''(v_i)=-\sum\limits_{j=1}^na_{u_j}^{(f)}v_j+a_{h}^{(v_i)}v_i-\\
           &-\sum\limits_{ 1\leq k<i}a_{s_{k,i}}^{(v_i)}v_k+\sum\limits_{ i<l\leq n}a_{s_{i,l}}^{(v_i)}v_l-\frac{\lambda^{(v_i)}}{2}v_i.
\end{split}
\end{equation}

Comparing the coefficients at the basis elements $f,\ z$ and $v_i$, (\ref{f4}) and (\ref{f5}) we obtain that

\begin{equation}\label{f6}
a_{u_i}^{(f)}=a_{h}^{(v_i)}-\frac{\lambda^{(v_i)}}{2}.
\end{equation}

Comparing (\ref{f3}) and (\ref{f6}) we obtain that
$$a_{u_i}^{(f)}=0,\quad a^{(v_i)}_{h}=\frac{\lambda_{z}^{(v_i)}}{2}.$$
So, $\Delta''(f)=0$ follows from equailty (\ref{asos}). We have the following connection $$a_{h}^{(v_i)}-\frac{\lambda_{z}^{(v_i)}}{2}=0$$
between the coefficients.
\end{proof}
\begin{lemma}\label{Lem6}
$\Delta''\left(e\right)=0$ and $a_{h}^{(u_i)}+\frac{\lambda^{(u_i)}}{2}=0$ in the formula (\ref{asos}).
\end{lemma}

\begin{proof}
From the equalities, take an element $x=e+\frac{1}{2}z+u_i$
\begin{equation}\begin{split}\label{f7}
\Delta''(x)&=[b,e+\frac{1}{2}z+u_i]+\mu\sigma(e+\frac{1}{2}z+u_i)=\\
             &=2b_{h}e-b_{f}h-\sum\limits_{j=1}^nb_{v_j}u_j+\frac{\mu}{2}z+\\
             &+b_{f}v_i+b_{h}u_i-b_{v_i}z+\sum\limits_{ 1\leq k<i}b_{s_{k,i}}u_k-\sum\limits_{ i<l\leq n}b_{s_{i,l}}u_l+\frac{\mu}{2}u_i.
\end{split}
\end{equation}

On the other hand,
\begin{equation}\begin{split}\label{f8}
\Delta''(x)&=\Delta''(e)+\Delta''(\frac{1}{2}z)+\Delta''(u_i)=-\sum\limits_{j=1}^na_{v_j}^{(e)}u_j+a_{h}^{(u_i)}u_i+\\
&+\sum\limits_{ 1\leq k<i}a_{s_{k,i}}^{(u_i)}u_k-\sum\limits_{ i<l\leq n}a_{s_{i,l}}^{(u_i)}u_l+\frac{\lambda^{(u_i)}}{2}u_i.
\end{split}
\end{equation}

Comparing the coefficients at the basis elements $e,\ z$ and $u_i$, (\ref{f7}) and (\ref{f8}) we obtain that

\begin{equation}\label{f9}
a_{v_i}^{(e)}=a_{h}^{(u_i)}+\frac{\lambda^{(u_i)}}{2}.
\end{equation}

From the equalities, take an element $x=e+\frac{1}{2}z-u_i$
\begin{equation}\begin{split}\label{f10}
\Delta''(x)&=[b,e+\frac{1}{2}z-u_i]+\mu\sigma(e+\frac{1}{2}z-u_i)=2b_{h}e-b_{f}h-\sum\limits_{j=1}^nb_{v_j}u_j+\frac{\mu}{2}z-\\
             &-b_{f}v_i-b_{h}u_i+b_{v_i}z-\sum\limits_{ 1\leq k<i}b_{s_{k,i}}u_k+\sum\limits_{ i<l\leq n}b_{s_{i,l}}u_l-\frac{\mu}{2}u_i.
\end{split}
\end{equation}

On the other hand,
\begin{equation}\begin{split}\label{f11}
\Delta''(x)&=\Delta''(e)+\Delta''(\frac{z}{2})-\Delta''(u_i)=-\sum\limits_{j=1}^na_{v_j}^{(e)}u_j-a_{h}^{(u_i)}u_i-\\
&-\sum\limits_{ 1\leq k<i}a_{s_{k,i}}^{(u_i)}u_k+\sum\limits_{ i<l\leq n}a_{s_{i,l}}^{(u_i)}u_l-\frac{\lambda^{(u_i)}}{2}u_i.
\end{split}
\end{equation}

Comparing the coefficients at the basis elements $e,\ z$ and $u_i$, (\ref{f10}) and (\ref{f11}) we obtain that
\begin{equation}\label{f12}
a_{v_i}^{(e)}=a_{h}^{(u_i)}+\frac{\lambda^{(u_i)}}{2}.
\end{equation}

Comparing (\ref{f9}) and (\ref{f12}) we obtain that
$$a_{v_i}^{(e)}=0,\quad a_{h}^{(u_i)}=-\frac{\lambda^{(u_i)}}{2}.$$
So, $\Delta''(e)=0$ follows from equailty (\ref{asos}). We have the following connection $$a_{h}^{(u_i)}+\frac{\lambda_{z}^{(u_i)}}{2}=0$$
between the coefficients.
\end{proof}

\begin{lemma}\label{Lem7}
$\Delta''\left(\mathfrak{so}_n\right)=\{0\}$ and  $\Delta''\left(\mathfrak{h}_n\right)=\{0\}.$
\end{lemma}

\begin{proof} Let $k,l(k\neq l)$ fixed number in set $\{1,2,...,n\}$.

From the equalities, take an element $x=u_k+s_{k,l}$, (if $l<k$ then $s_{k,l}=-s_{l,k}$)
\begin{equation}\begin{split}\label{f13}
\Delta''(x)&=[b,u_k+s_{k,l}]+\mu\sigma(u_k+s_{k,l})=\\
             &=b_{f}v_k+b_{h}u_k-b_{v_k}z+\sum\limits_{ 1\leq p<k}b_{s_{p,k}}u_p-\\
             &-\sum\limits_{ k<q\leq n}b_{s_{k,q}}u_q+\frac{\mu}{2}u_k-b_{u_l}u_k+b_{u_k}u_l-b_{v_l}v_k+b_{v_k}v_l+\\
             &+\sum\limits_{ 1\leq p< k}b_{s_{p,k}}s_{p,l}+\sum\limits_{ l<q\leq n}b_{s_{l,q}}s_{q,k}+
             \sum\limits_{ 1\leq p< l}b_{s_{p,l}}s_{p,k}+\sum\limits_{ k<q\leq n}b_{s_{k,q}}s_{l,q}.
\end{split}
\end{equation}

On the other hand,
\begin{equation}\begin{split}\label{f14}
\Delta''(x)&=\Delta''(u_k)+\Delta(s_{k,l})=\sum\limits_{ 1\leq j<k}a_{s_{j,k}}^{(u_k)}u_j-\sum\limits_{ k<j\leq n}a_{s_{k,j}}^{(u_i)}u_j-\\
&-a_{u_l}^{(s_{k,l})}u_k+a_{u_k}^{(s_{k,l})}u_l-a_{v_l}^{(s_{k,l})}v_k+a_{v_k}^{(s_{k,l})}v_l.
\end{split}
\end{equation}

Comparing the coefficients at the basis elements $ z$ and $y_l$, (\ref{f13}) and (\ref{f14}) we obtain that
\begin{equation*}\begin{split}\label{f15}
a_{v_k}^{(s_{k,l})}=0.
\end{split}
\end{equation*}

Similarly, from equality:
\begin{equation*}\begin{split}\label{f16}
&\Delta''\left(u_l+s_{l,k}\right)=\Delta''\left(u_k\right)+\Delta''\left(s_{l,k}\right) \text{ we obtain }a_{v_l}^{(s_{k,l})}=0;\\
&\Delta''\left(v_k+s_{k,l}\right)=\Delta''\left(v_k\right)+\Delta''\left(s_{k,l}\right) \text{ we obtain } a_{u_k}^{(s_{k,l})}=0;\\
&\Delta''\left(v_l+s_{l,k}\right)=\Delta''\left(v_k\right)+\Delta''\left(s_{l,k}\right) \text{ we obtain }a_{u_l}^{(s_{k,l})}=0.\\
\end{split}
\end{equation*}

Comparing the coefficients at the basis elements $x_j (j\neq k, j \neq l) $ and $s_{j,l}$, (\ref{f13}) and (\ref{f14}) we obtain that
\begin{equation}\begin{split}\label{f17}
a_{s_{k,j}}^{(u_k)}=0, (j \neq l).
\end{split}
\end{equation}

Take a number $i$ which $i\neq k$ and $i\neq l$.
Similarly, from equality:
\begin{equation}\begin{split}\label{f18}
&\Delta''\left(u_k+s_{k,i}\right)=\Delta''\left(u_k\right)+\Delta''\left(s_{k,i}\right) \text{ we obtain } a_{s_{k,l}}^{(u_k)}=0;\\
&\Delta''\left(v_k+s_{k,l}\right)=\Delta''\left(v_k\right)+\Delta''\left(s_{k,l}\right) \text{ we obtain } a_{s_{k,j}}^{(v_k)}=0, (j \neq l);\\
&\Delta''\left(v_k+s_{k,i}\right)=\Delta''\left(v_k\right)+\Delta''\left(s_{k,i}\right) \text{ we obtain } a_{s_{k,l}}^{(v_k)}=0.\\
\end{split}
\end{equation}

Thus, according to  (\ref{asos}), (\ref{f17}) and (\ref{f18}) we have
$$\Delta''\left(\mathfrak{so}_n\right)=\{0\}\quad \text{and}\quad  \Delta''\left(\mathfrak{h}_n\right)=\{0\}.$$
\end{proof}

Now we are in position to prove Theorem \ref{thm21}.

\textit{Proof of Theorem \ref{thm21}:}
From \eqref{L11} and lemmas (\ref{Lem5} -- \ref{Lem7}) we obtain
\begin{equation}\label{f19}
\Delta''=0.
\end{equation}
Together \eqref{L9} and \eqref{f19} gives
\begin{equation}\label{f20}
\Delta'=\mathrm{ad}(x_0)+\lambda^{(h)}\sigma.
\end{equation}
Together \eqref{LDer1} and \eqref{f20} gives
\begin{equation*}\label{f21}
\Delta=D+\mathrm{ad}(x_0)+\lambda^{(h)}\sigma.
\end{equation*}
Hence, any local derivation of the algebra $\ms_n(n\ge 3)$ is a derivation.
\hfill $\Box$

\end{document}